\documentclass[11pt]{article}

\usepackage{makeidx}
\usepackage{sec}
\usepackage{pb-diagram}
\usepackage{lamsarrow}
\usepackage{pb-lams}
\usepackage{amsmath}
\usepackage{amsthm}
\usepackage{epsfig}
\usepackage{amssymb}

\theoremstyle{plain}
\newtheorem{theorem}{Theorem}[section]
\newtheorem{lemma}[theorem]{Lemma}
\newtheorem{proposition}[theorem]{Proposition}
\newtheorem{corollary}[theorem]{Corollary}
\newtheorem{definition}[theorem]{Definition}

\theoremstyle{definition}

\newtheorem{remark}[theorem]{Remark}

\DeclareMathOperator{\Hom}{{\rm Hom}}

\newcommand\bC{{\mathbb C}}

\newcommand{\bQ}{{{\mathbb Q}}}
\newcommand\bR{{\mathbb R}}
\newcommand\bS{{\mathbb S}}

\newcommand\bZ{{\mathbb Z}}

\renewcommand{\Box}{\blacksquare}

\newcommand{\gq}{{\mathfrak q}}

\newcommand{\gR}{{\mathfrak R}}

\newcommand{\X}{{\mathfrak X}}

\newcommand{\sw}{{\bf S}{\bf W}}

\newcommand\aug{\mathfrak{aug}}

\newcommand{\dir}{{\mathfrak D}}

\newcommand{\A}{\mathcal A}

\newcommand{\f}{\mathcal F}

\renewcommand{\k}{\mathcal K}

\newcommand{\p}{{\mathcal P}}
\newcommand{\q}{\mathcal Q}
\newcommand{\R}{\mathcal R}

\renewcommand{\t}{{\mathcal T}}

\newcommand{\ra}{\rightarrow}
\newcommand{\hra}{\hookrightarrow}
\newcommand{\Lra}{{\longrightarrow}}

\newcommand{\lan}{\langle}
\newcommand{\ran}{\rangle}

\def\inpr{\mathbin{\hbox to 6pt{\vrule height0.4pt width5pt depth0pt \kern-.4pt \vrule height6pt width0.4pt depth0pt\hss}}}

\newcommand{\si}{{\sigma}}

\newcommand{\ve}{{\varepsilon}}
\newcommand{\vfi}{{\varphi}}

\newcommand{\pa}{\partial}

\newcommand{\eps}{\epsilon}
\newcommand{\ssw}{{\bf sw}}

\makeindex

\begin{document}

\title{Seiberg-Witten invariants of rational homology $3$-spheres}

\author{Liviu I. Nicolaescu\thanks{This work was partially supported by NSF grant DMS-0071820.} \\University of Notre Dame\\Notre Dame, IN 46556\\http://www.nd.edu/ $\tilde{}$ lnicolae/ }

\date{Version 1, February 2001}

\maketitle

\begin{abstract} We prove that the Seiberg-Witten invariants of a rational homology sphere are determined by the Casson-Walker invariant and the Reidemeister torsion.
\end{abstract}

\section*{Introduction}
 In 1996 Meng and Taubes \cite{MT} have established    a  relationship between the Seiberg-Witten  invariants of a (closed) $3$-manifold with $b_1\geq 0$.    A bit later    Turaev \cite{Tu6,Tu}     refined    Meng-Taubes' result, essentially identifying these two invariants for such $3$-manifolds. In \cite{Tu6} Turaev left open the question  of establishing a connection  between these two invariants in the remaining case, that of rational homology spheres.

Around the same time, Lim \cite{Lim2} succeeded in providing a combinatorial description of the Seiberg-Witten invariants of integral homology spheres. Namely, in this case they  coincide with the Casson invariant.    In \cite{Nico3} we  investigated a special class of rational homology spheres, the lens spaces, and we proved  that the Seiberg-Witten invariants of such spaces  are equivalent to the Casson-Walker invariant and the Reidemeister torsion, and raised the question whether this is the case in general.    Recently Marcolli  and Wang \cite{MW} (see also the  related work of Ozsv\'{a}th-Zsab\'{o} \cite{OS}) have shown that the Seiberg-Witten invariants     of a ${\bQ}HS$ determine the Casson-Walker invariant. Additionally, they have proved  a very general surgery formula involving the Seiberg-Witten invariants.

In the present paper we prove that for rational homology spheres we have
\[
\mbox{SW}\Longleftrightarrow \mbox{modified Reidemeister torsion} \stackrel{def}{:=}\mbox{Casson-Walker $+$ Reidemeister torsion}.
\]
Our strategy is based on surgery formul{\ae} for Seiberg-Witten invariants developed   in \cite{MW}, and surgery formul{\ae} for the modified torsion,   described in \cite{Tu, Walker}.

Both the Seiberg-Witten invariant and the modified Reidemeister torsion can be thought  of as  ${\bQ}$-valued functions
on the first homology group $H$ of  a given rational homology sphere $M$.  We denote by $D_M$ the difference    of these two functions.      Proving the equality of these two invariants is equivalent to showing that $D_M\equiv 0$.

We found it extremely convenient to work not with $D_M$ but with
its Fourier transform $\widehat{D}_M: H^\sharp\ra {\bC}$, where
$H^\sharp$ is the dual of $H$.  For example, Marcolli-Wang result
\cite{MW} translates into $\widehat{D}_M(1)=0$, for all  rational
homology spheres. Additionally, the true nature of the  surgery
formul{\ae} is better displayed in the Fourier picture. To explain
the gist  of these formul{\ae} consider a $3$-manifold $N$ with
$b_1=1$ and boundary $T^2$. $N$ can be thought of as the
complement of a knot in a ${\bQ}HS$.  Pick two simple  closed
curves  $c_1$, $c_2$ on $\pa N$ with nontrivial intersection
numbers with the longitude $\lambda\in H_1(\pa N, {\bZ})$.

By Dehn surgery with $c_i$ as attaching curves we obtain two rational
homology spheres  $M_1$, $M_2$ and two knots $K_i\hra M_i$,
$i=0,1$.    Let $H_i:= H_1(M_i, {\bZ})$, $G:=H_1(N,\pa N; {\bZ})$. The knot $K_i$ determines
a subgroup $K_i^\perp \subset H_i^\sharp$, the characters vanishing on $K_i$.   These subgroups are naturally isomorphic to $G$  and thus we have a natural isomorphism
\[
f: K_1^\perp\ra K_2^\perp.
\]
The surgery formul{\ae}   have the form
\[
\lan \lambda, c_2\ran\widehat{D}_M(f^*\chi)=\lan \lambda, c_1\ran \widehat{D}_{M_2}(\chi)+ |G| \k,\;\;\forall \chi\in K_2^\perp
\]
where $\lan\bullet,\bullet\ran$ denotes the intersection pairing on $H_1(\pa N, {\bZ})$, and $\k$ is an universal correction term which depends only on  the divisibility $m_0$  of  the longitude, and the  $SL_2({\bZ})$-orbit  of the pair $(c_1,c_2)$ with respect to the  obvious action of this group on the space of pairs of primitive vectors in a $2$-dimensional lattice. We will thus write $\k_{m_0; [c_1,c_2]}$, and call the triplet $(m_0;[c_1,c_2])$ the arithmetic type of the surgery. The results of \cite{OS} prove that
\[
\k_{1;[c_1,c_2]}\equiv 0,\;\;\forall[c_1,c_2].
\]
We call surgeries with $m_0=1$ primitive. The admissible surgeries have trivial correction
term.  We denote  by $\X$ the class of rational  homology spheres
$M$ such that $\widehat{D}_M\equiv 0$.  Both the family of
admissible surgeries and the family $\X$ are ``time dependent''
families, and  during our proof we will gradually produce larger
and larger classes of surgeries/ manifolds   inside these
families.

The class $\X$ is  closed   under connected sums and  certain primitive surgeries (see \ref{ss: 41}). Using this preliminary  information we are able to show    that all  homology lens spaces belong to $\X$. The proof uses   Kirby calculus, and  we learned  it from Nikolai Saveliev. As a bonus, we can include many more  arithmetic types  of Dehn surgeries in the class of admissible surgeries.

Loosely speaking, the  homology lens spaces have the simplest
linking forms.  We take this idea seriously, and we define an
appropriate notion of complexity of a linking form.    The proof
then proceeds by induction, including in $\X$ manifolds of larger
and larger complexity. This process  also increases the class of
admissible surgeries,  which can be used at the various inductive
steps.        Such a proof is feasible if we can produce   a large
supply of complexity reducing   Dehn surgeries.   Fortunately, this is   can be seen  using  elementary arithmetic.

Our proof    also shows  that the  invariant  introduced  by Ozsv\'{a}th and Szab\'{o} in \cite{OS} also coincides with  Casson-Walker $+$ Reidemeister torsion, thus answering a question raised in that paper.      Moreover    our result  establishes connections between the Kreck-Stolz invariant \cite{KS}, the linking form, and  an invariant we introduced  in  \cite{N3}.

\bigskip

\noindent {\bf Acknowledgments}     I  learned  about the
Reidemeister torsion from  Frank Connolly. He  explained  the
basics to me,           and  over countless hours of discussions
he  listened and answered  patiently to    many  of my questions and   my
half-developed ideas. These discussions    always enriched my knowledge, and made me aware   of many of the subtleties of this subject.

I am grateful to Nikolai Saveliev for his patience and expertise in answering   the questions I had while working on  this project. In particular, his proof of Lemma \ref{lemma: sav}    got  me  ``out of a hole''.  I want to thank my colleagues Stephan Stolz and Larry Taylor for  many illuminating discussions.

\bigskip

\noindent {\bf Basic Notations and Terminology } A closed, compact, oriented $3$-manifold will be denoted by $M$. We will set $H=H_1(M, {\bZ})\cong H^2(M, {\bZ})$, and we will  denote the group operation multiplicatively. We denote by $T(H)$ the torsion part of $H$ and by $T_2(H)$ the 2-torsion part
 \[
 T_2(H):=\bigl\{ h\in T(H);\;\;h^2=1\bigr\}.
 \]
We set $\Theta=\Theta_M:=\sum_{h\in T(H)}h\in {\bZ}[H]$. For any $P=\sum_{h\in H}P_hh\in {\bZ}[H]$ we set
\[
\bar{P}=\sum_{h\in H}P_hh^{-1}.
\]
The letter $N$ will be reserved for  compact, oriented three manifolds with boundary $\pa N\cong T^2$ such that $b_1(N)=1$. Equivalently, $N$ can be viewed as  the complement of a knot in a rational homology sphere. We set $G=H_1(N,\pa N)\cong H^2(N, {\bZ})$.

We will denote   by $Spin^c(M)$ the space  of isomorphism classes of $spin^c$ structures on $M$.  We will denote a generic $spin^c$ structure by $\si$. $Spin^c(M)$ is a $H$-torsor, and we will denote by
\[
Spin^c(M)\times H\ni(\si, h)\mapsto h\si
\]
the natural  action of $H$ on $Spin^c(M)$.   The natural involution on $Spin^c(M)$ will be denoted by $\si\mapsto\bar{\si}$.  The complex line bundle associated to $\si$ will be denoted by $\det(\si)$.  We can identify $\det(\si)$ via the first Chern class  with an element in $H$. Note that
\[
\det(h\si)=h^2\det(\si).
\]
We will denote by $Spin(M)$  the space of  isomorphism classes of $spin$-structures on $M$. A generic $spin$ structure   will be denoted by $\epsilon$. $Spin(M)$ is naturally a  $T_2(H)$-torsor.  We use the same notation to denote the action of $T_2(H)$ on $Spin(M)$. Every $spin$ structure $\epsilon$ induces a canonical  $spin^c$-structure  $\si(\epsilon)$. Moreover
\[
\si(\epsilon)=\overline{\si(\epsilon)},\;\;\si(h\epsilon)=h\si(\epsilon),\;\;\forall h\in T_2(H).
\]
For any Abelian group $A$ we will denote by $A^\sharp$ its dual, $A^\sharp=\Hom(A, S^1)$. Finally for every $\chi\in H^\sharp$, and any $P=\sum_{h\in H}P_hh\in {\bC}[H]$ we set
\[
\hat{P}(\chi):=\sum_{h\in H}P_h\chi(h)\in {\bC}.
\]
The  function $H^\sharp\ni \chi\mapsto\hat{P}(\chi)$ is essentially the Fourier transform of $P$. Note that
\[
\hat{P}(1):=\sum_{h\in H}P_h.
\]
Moreover
\[
\hat{\Theta}_M(1)=|T(H)|,\;\;\hat{\Theta}_M(\chi)=0,\;\;{\rm if}\;\;\chi\neq 1, \;\;{\rm and}\;\;\exists m>1\;\;\chi(H)\subset U_m.
\]

 For every positive integer $m$ we denote by $U_m\subset S^1$ the multiplicative group of $m$-th roots of $1$.


\tableofcontents

\section{The modified Seiberg-Witten invariants of $3$-manifolds}
We  want to present in a form suitable for our goals, some basic structural facts concerning the Seiberg-Witten invariants of $3$-manifolds. For more details we refer to \cite{CMW, MT, Lim}.

\subsection{The case $b_1>1$.}
Suppose $b_1(M)>1$.   Fix an orientation on $H\otimes {\bR}$.  The Seiberg-Witten invariant   of $M$ is a function
\[
\ssw_M: Spin^c(M)\ra {\bZ}.
\]
$\ssw_M(\si)$ is a signed count of $\si$-monopoles,   objects determined by additional geometric data on $M$,  Riemann metric $g$ and a closed $2$-form $\eta$.   The   chosen orientation on $H\otimes {\bR}$   associates a sign to each monopole, and the signed count is independent of $g$ and $\eta$. The Seiberg-Witten invariant has the following properties.

\medskip

\noindent $\bullet$ $\ssw_M(\si)=0$ for all but finitely many $\si$'s.

\noindent $\bullet$ $\ssw_M(\si)=\ssw_M(\bar{\si})$, $\forall \si$.

\medskip

For every $\si$ we can form the element
\[
\sw_{M,\si}\in{\bZ}[H],\;\;\sw_{M,\si}=\sum_{h\in H}\ssw_M(h^{-1}\si)h.
\]
Note  that for every $h_0\in H$ we have
\[
\sw_{M,h_0\si}=h_0\sw_{M,\si}.
\]
Moreover
\[
\sw_{M,\si}=\det(\si)\sw_{M,\bar{\si}}=\det(\si)\overline{\sw}_{M,\si}.
\]
In particular, for any $spin$ structure $\eps$ we have $\det(\si(\eps))=1$ so that
\[
\sw_{M,\si(\eps)}= \overline{\sw}_{M,\si(\eps)}.
\]
For simplicity we set $\sw_{M,\eps}:=\sw_{M,\si(\eps)}$, $\forall\eps$.

\subsection{The case $b_1=1$.}
Suppose $b_1(M)=1$, and fix an orientation on $H\otimes {\bR}$. In this case, a choice of orientation is determined by fixing an isomorphism $H\otimes {\bR}\ra {\bR}$.  To  describe the Seiberg-Witten invariant of $M$ we need to recall the rudiments of its construction.
theory. Fix a metric  $g$. The chosen orientation on $H\otimes {\bR}$ defines  a harmonic $1$-form $\omega_g$   such that $\omega_g$ induces the  chosen orientation on $H\otimes {\bR}$, and $\|\omega_g\|_{L^2(g)}=1$.   By rescaling the metric  we can assume that  $\omega_g$ also generates the image of $H$ in $H\otimes {\bR}$.

Note  that  the chosen orientation  produces an isomorphism $H/(T(H)\ra {\bZ}$, and thus a map
\[
\deg : H\ra H/T(H)={\bZ}.
\]

For $\si\in Spin^c(M)$ denote by $\p_\si(g)$ the space of closed $2$-forms such that
\[
w(\si,\eta):=\int_M\omega_g\wedge \eta-2\pi c_1(\det\si) \neq 0.
\]
It is decomposed into two chambers
\[
\p_\si^\pm(g) =\Bigl\{ \eta\in \p_\si(g);\;\;\pm w(\si,\eta)>0 \Bigr\}.
\]
For $\eta\in \p_\si(g)$ we denote by $\ssw_M^\pm(\si,\eta)$ the signed count of $(\si,g,\eta)$-monopoles.  It is known that
\[
\ssw_M(\si,\eta)=\ssw_M(\bar{\si},\eta)
\]
$\ssw_M(\si,\eta)=0$ for all but finitely many $\si$'s, and
\[
\ssw_M(\si,\eta_1)=\ssw_M(\si,\eta_2),\;\;{\rm if}\;\;w(\si, \eta_1)\cdot w(\si,\eta_2)>0.
\]
We set
\[
\ssw_M^\pm (\si):= \ssw_M(\si,\eta),\;\;\pm w(\si,\eta)>0.
\]
The wall crossing formula (see \cite{Lim}) states that
\[
\ssw_M^+(\si) -\ssw_M^-(\si) =\frac{1}{2}\deg(\det\si).
\]
Set
\[
\sw_{M,\si,\eta}=\sum_{h\in H}\ssw_M(\si,\eta)\in {\bZ}[H],\;\;\sw_{M,\si}=\sum_{h\in H}\ssw_M^+(h^{-1}\si)h\in {\bZ}[[H]].
\]
Suppose we  pick $\si=\si(\eps)$ and $\eta=\eta_0$ such that $\int_M\omega\wedge \eta_0$ is a very small positive number. Fix $T\in H$ such that $\deg(T)=1$. We can rephrase the wall crossing formula in the more compact form
\[
\sw_{M,\si(\eps)}=\sw_{M,\si(\eps),\eta_0}+ \frac{\Theta_MT}{(1-T)^2}.
\]
We set $W_M:=\frac{\Theta_MT}{(1-T)^2}\in {\bZ}[[H]]$, and we will refer it as {\em  wall-crossing correction}\footnote{Observe that the wall-crossing correction  $\Theta_MT(1-T)^{-2}\in {\bZ}[[H]]$ is independent of the choice of $T$ in $H$ such that $\deg T=1$.} of $M$. Observe that
\[
\sw_{M,\si(\eps)}^0:=\sw_{M,\si(\ve)}-W_M=\sw_{M,\si(\eps),\eta_0}\in {\bZ}[H]
\]
satisfies the  symmetry condition
\[
\sw^0_{M,\si(\eps)}=\overline{\sw}^0_{M,\si(\eps)},
\]
and
\[
\sw^0_{M,\si(h_0\eps)}=h_0\sw^0_{M,\si(\eps)},\;\;\forall h_0\in T_2(H).
\]
We will refer to $\sw^0_{M,\si(\eps)}$ as {\em the modified  Seiberg-Witten invariant of $M$.}

\subsection{The case $b_1=0$.}
Suppose now  that $b_1(M)=0$, i.e. $M$ is a rational homology sphere.  Fix $\si\in  Spin^c(M)$.  In this case  the signed count of $(\si, g,\eta)$-monopoles depends on $(g,\eta$ in a more complicated way.  To  produce a topological invariant we need to  add a correction to this count.  For simplicity, we describe this correction only  when $\eta=0$.

Denote by ${\bS}_\si$ the  bundle of complex spinors determined by $\si$. The line bundle $\det\si=\det{\bS}_\si$ admits an unique equivalence class of flat connections.  Pick one such flat connection $A_\si$ and denote by $\dir_{A_\si}$   Dirac operator  on ${\bS}_\si$  determined by the twisting connection $\si$.   We denote its eta invariant  by $\eta_{dir}(g,\si)$. Also, denote by $\eta_{sign}(g)$ the eta invariant of the odd signature operator  determined by $g$. Finally define the Kreck-Stolz invariant of $(g,\si)$ by
\[
KS(g,\si)=4\eta_{dir}(g,\si)+\eta_{sign}(g).
\]
Define {\em the modified Seiberg-Witten invariant} of $(M,\si)$ by
\[
\ssw_M^0(\si)= \frac{1}{8}KS(g,\si)+\ssw_M(\si)\in {\bQ}.
\]
As shown in \cite{Lim}, the above quantity is  independent of the metric, and it is a topological invariant. Set
\[
\sw^0_{M,\si}:=\sum_{h\in H}\ssw_M(h^{-1}\si)h\in {\bQ}[H].
\]
If $\si=\si(\eps)$ we have
\[
\sw^0_{M,\si(\eps)}=\overline{\sw}^0_{M,\si(\eps)}.
\]

\subsection{Summary}
\label{ss: 14}
Let us  coherently organize the facts    explained so far.  We say that $M$ is {\em  homologically oriented} if $H\otimes {\bR}$ is oriented. The modified  Seiberg-Witten invariant   associates to each closed, compact, homologically  oriented $3$-manifold $M$, and each $\eps\in Spin(M)$  a  ``Laurent polynomial'' $\sw^0_{M,\eps}\in {\bQ}[H]$ with the following properties.
\begin{equation}
\sw^0_{M,\eps}\in {\bZ}[H],\;\;{\rm if}\;\;b_1(M)>0,
\label{eq: sw1}
\end{equation}
\begin{equation}
\sw^0_{M,\eps}=\overline{\sw}^0_{M,\eps},
\label{eq: sw2}
\end{equation}
and
\begin{equation}
\sw^0_{M,h_0\eps}=h_0\sw^0_{M,\eps},\;\;\forall h_0\in T_2(H).
\label{eq: sw3}
\end{equation}

\section{The modified Reidemeister-Turaev  torsion of $3$-manifolds}
\setcounter{equation}{0}
In this section we survey   the results of V. Turaev \cite{Tu2,Tu4,Tu5,Tu6, Tu} in a language appropriate to our goals.

\subsection{Turaev's refined torsion}
 The Reidemeister-Turaev torsion associates to each   of a homologically oriented $3$-manifold  $M$, and  each $spin^c$ structure $\si$ on $M$ a ``formal power series'' $\t_{M,\si}\in {\bQ}[[H]]$ with the following properties.
 \[
 \t_{M,\si}\in {\bZ}[[H]],\;\;b_1(M)>1
 \]
 \[
 (1-T)^2\t_{M,\si}\in {\bZ}[H],\;\;b_1(M)=1,\;\;\deg T=1
 \]
 \[
 \t_{M,\si}\in {\bQ}[H],\;\;,\;\;\t_{M,\si}(1)=0,\;\;b_1(M)=0.
 \]
 Moreover
 \[
 \t_{M,h_0\si}=h_0\t_{M,\si},\;\;\forall h_0\in H,
 \]
 and
 \[
 \t_{M,\si}=\det(\si)\overline{\t}_{M,\si}.
 \]
 For $\eps\in Spin(M)$ set $\t_{M,\eps}=\t_{M,\si(\eps)}$. It follows that
 \[
 \t_{M,\eps}=\overline{\t}_{M,\eps}.
 \]
Using \cite{Tu5,Tu6} and \cite[Appendix 3]{Tu} we deduce that when $b_1(M)=1$ we have
\[
\t_{M,\eps}^0:=\t_{M,\eps}-W_M\in {\bZ}[H],
\]
and moreover
\[
\t^0_{M,\eps}=\overline{\t}^0_{M,\eps}.
\]
When $b_1(M)=0$ denote by $CW(M)$ the Casson-Walker invariant of $M$ and define
\[
\t^0_{M,\eps}=\t_{M,\eps}+\frac{1}{2}CW(M)\Theta_M.
\]
Observe that $\t^0_{M,\eps}(1)=\frac{1}{2}|H|CW(M)$=  Lescop invariant of $M$ (see \cite[p. 80]{Lescop}).

We we will refer to  the quantities $\t_{M,\eps}^0$ for $b_1(M)=0,1$ the {\em modified  Reidemeister-Turaev torsion of  $M$.} For uniformity, we set $\t^0_M=\t_M$ when $b_1(M)>1$. Summarizing, we conclude that the modified Reidemeister-Turaev torsion  associates to each homologically oriented $3$-manifold $M$, and to each $spin$ structure $\eps$ on $M$ a ``Laurent polynomial'' $\t^0_{M,\eps}\in {\bQ}[H]$ with the following properties.
\begin{equation}
\t^0_{M,\eps}\in {\bZ}[H],\;\;{\rm if}\;\;b_1(M)>0,
\label{eq: tor1}
\end{equation}
\begin{equation}
\t^0_{M,\eps}=\overline{\t}^0_{M,\eps},
\label{eq: tor2}
\end{equation}
and
\begin{equation}
\t^0_{M,h_0\eps}=h_0\t^0_{M,\eps},\;\;\forall h_0\in T_2(H).
\label{eq: tor3}
\end{equation}

\subsection{Relations between the torsion and the Seiberg-Witten invariant}
The Seiberg-Witten invariant and the modified Reidemeister torsion  are related. More precisely we have the following result.

\begin{theorem}(a)  $\sw^0_{M,\eps}=\t^0_{M,\eps}$ if $b_1(M)>0$; see {\bf \cite{MT,Tu6, Tu}}.

\noindent (b)$\sw^0_M(1)=\t^0_M(1)$ if $b_1(M)=0$; see  {\bf \cite{CMW, Lim2}}.

\noindent (c)$\sw^0_M=\t_M^0$ if $M$ is a lens space; see {\bf \cite{Nico3}}.
\label{th: relat}
\end{theorem}
Part (c) of the above theorem can be slightly strengthened to
\begin{equation}
\sw_M^0=\t_M^0,\;\;\mbox{if $M$ is a connected sum of lens spaces}.
\label{tosw}
\end{equation}
This equality follows from the  vanishing of the torsion  under connected sums, the additivity of the  Casson-Walker invariant, and the additivity of the Kreck-Stolz invariant. (This follows from the very general surgery results for eta invariants in \cite{KL}.)

Later on we will need the following consequence of  Theorem \ref{th: relat} (a).

\begin{proposition} If $M$ is a homologically oriented  $3$-manifold such that $b_1(M)=1$ then
\[
\widehat{\t}^0_M(1)=\widehat{\sw}^0_M(1)=\frac{1}{2}\Delta_M''(1),
\]
where $\Delta_M\in {\bZ}[[T^{1/2}, T^{-1/2}]]$ denotes the symmetrized Alexander polynomial of $M$ normalized such that $\Delta_M(1)=|T(H)|$.
\label{prop: alex}
\end{proposition}

\noindent {\bf Proof}\hspace{.3cm} The  projection $H\ra H/T(H)={\bZ}$ induces  a morphism
\[
\aug: {\bZ}[[H]]\ra {\bZ}[t,t^{-1}]]
\]
called augmentation.  Fix $T\in H$ such that $\deg T=1$. The  symmetrized Alexander polynomial $\Delta_M$   is uniquely determined by the condition
\[
\aug\t_{M,\eps}= T^{k/2}\frac{\Delta_M(T)}{(1-T)^2}
\]
for some $k\in {\bZ}$.  Using Theorem \ref{th: relat}(a) we deduce
\[
T^{k/2}\frac{\Delta_M(T)}{(1-T)^2}=\aug\sw_M=\aug\sw_M^0 +\aug(\Theta_M)\frac{T}{(1-T)^2}=\aug\sw_M^0 +|H|\frac{T}{(1-T)^2}.
\]
We conclude that
\[
T^{k/2 -1}\Delta_M(T)= (T-2+T^{-1})\aug\sw^0_M(T)+|H|.
\]
The symmetry of $\sw^0$ implies $\sw^0_M(T)=\sw^0_M(T^{-1})$, and since $\Delta_M$ satisfies a similar  symmetry we conclude $k/2-1=0$. Hence
\[
\Delta_M(T)=(T-2+T^{-1})\aug\sw^0_M(T)+|H|.
\]
Differentiating the above equality twice at $T=1$ we deduce
\[
\Delta_M''(1)= 2\widehat{\aug\,\sw}_M(1)=2\widehat{\sw}^0(1).\;\;\;\Box
\]

\begin{remark}  Observe a  nice ``accident''. Suppose $M$ is as in Proposition \ref{prop: alex}. Then
\[
W_M=\Theta_M\sum_{n\geq 1}nT^n.
\]
Formally
\[
\widehat{W}_M(1)= \widehat{\Theta}_M(1)\sum_{n\geq 1}n =|H|\sum_{n\geq 1}n=|H|\zeta(-1)=-\frac{1}{12}|H|,
\]
where $\zeta(s)$ denotes Riemann's zeta function.  In particular
\[
\widehat{\sw}_M(1)=\widehat{\sw}_M^0(1)+\widehat{W}_M(1)=\frac{1}{2}\Delta_M''(1)-\frac{1}{12}|H|.
\]
The expression in the right-hand-side  is precisely the Lescop invariant of $M$.
\end{remark}

We can now state the main result of this paper.

\begin{theorem}
\[
\sw_M^0=\t_M^0
\]
for any, closed, compact, oriented $3$-manifold $M$.
\label{th: main}
\end{theorem}

\section{Surgery formul{\ae}}
\setcounter{equation}{0}

\subsection{Dehn surgery}
We want to survey a few basic facts concerning Dehn surgery. For more details and examples we refer to \cite{N2}.

Consider a $3$-manifold $N$ as in the introduction, i.e. $b_1(N)=1$, $\pa N\cong T^2$, and set $G:=H_1(N,\pa N; {\bZ})$. We orient $\pa N$ as boundary of $N$ using the  outer-normal first convention.  Denote by ${\bf j}$ the inclusion induced morphism
\[
{\bf j}: H_1(\pa N, {\bZ})\ra H_1(N, {\bZ}).
\]
The kernel of ${\bf j}$  is a rank one Abelian group. We can select a generator $\lambda$ of $\ker {\bf j}$ by specifying an orientation on $H^1(N, {\bZ})\cong H_2(N, \pa N; {\bZ})$.  We can write $\lambda=m_0\lambda_0$ where $m_0>0$ and $\lambda_0\in H_1(\pa N,{\bZ})$ is a primitive class.   $\lambda$ is called the longitude of $N$ and $m_0$ is called  the {\em divisibility} of $N$.Fix $\mu_0\in H_1(\pa N, {\bZ})$ such that $\lambda_0\cdot \mu_0=1$, where the dot denotes the intersection pairing on $H_1(\pa N, {\bZ})$.

Denote by $X$ the solid torus  $S^1\times D^2$, so that $\pa X=T^2$ Set $\ell_0= S^1\times \{ {\bf pt}\}$ and $ {\bf m}_0=\{ {\bf pt}\}\times \pa D^2$. We regard $ \ell_0$ and ${\bf m}_0$ as elements in $H_1(\pa X, {\bZ})$. They satisfy ${\bf m}_0\cdot \ell_0=1$. Fix an orientation reversing   diffeomorphism  $\Gamma:\pa X\ra \pa N$ such that
\[
\Gamma_*({\bf m}_0)=\mu_0,\;\;\Gamma_*(\ell_0)=\lambda_0.
\]
Every    $\vfi\in SL_2({\bZ})$  determines an isotopy class of orientation preserving  diffeomorphisms of $T^2$. We can for a closed $3$-manifold
\[
M_\vfi:= X\coprod_{\Gamma\circ\vfi:\pa X\ra \pa N}N.
\]
We say that $M_\vfi$ is obtained  by {\em Dehn surgery} with gluing map
$\vfi$. The integer $m_0$ is called the divisibility of the surgery. The    manifold $M_\vfi$ is uniquely determined up to a
diffeomorphism by  the attaching  curve $c=\Gamma\circ\vfi({\bf
m}_0$. We can write $c=c_{p/q}:=p\mu_0+ q\lambda_0$, $(p,q)=1$. The
diffeomorphism type of $M_\vfi$ is   uniquely determined  by the
ration $p/q$.  Instead of $M_{\vfi}$ we will write $M_{p/q}$. We set $H_{p/q}:= H_1(M_{p/q}, {\bZ})$. The  core of the solid torus determines an element $K_{p/q}\in H_{p/q}$.

We want to point out that the integer $q$ depends on the choice of $\mu_0$ while $p$ is invariantly determined  by the equality
\[
p:= \lambda_0\cdot c.
\]
We  we refer to $p$ as the {\em multiplicity} of the surgery.

The  group $H_{p/q}$ is determined from the short exact sequence
\[
0\ra \lan {\bf j}c_{p/q}\ran \ra H_1(N, {\bZ})\ra H_{p/q}\ra 0.
\]
We also have canonical isomorphisms
\[
\Phi_{p/q}: G\ra H_{p/q}/\lan K_{p/q}\ran=: R_{p/q}.
\]
We obtain a natural  projection $\pi_{p/q}:H_{p/q}\ra G$.

The long exact sequence of the pair $(N,\pa N)$ implies
\[
G= H_1(N, {\bZ})/{\bf j} H_1(\pa N, {\bZ}).
\]
We deduce the following result.

\begin{lemma} The characters  of $G$   are precisely the characters of $H_1(N, {\bZ})$ which vanish on ${\bf j}H_1(\pa N, {\bZ})$.  Also, we can think of the characters of $G$ as characters $\chi$ of $H_{p/q}$  such that $\chi(K_{p/q})=1$.
\label{lemma: surg}
\end{lemma}

When $p\neq 0$, $H_{p/q}$ is a finite Abelian group and
\[
|H_{p/q}|= pm_0 |G|.
\]
In this case, we denote by ${\bf lk}_{p/q}$ the linking form of $M_{p/q}$.

Observe that  $b_1(M_{0/1})=1$.  $K_{0/1}$ can be written   as $m_0h$ where $h\in H_{0/1}$ generates the free part of $H_{0/1}$. $M_{0/1}$ carries a natural homology orientation, induced from the orientation of $H^1(N, {\bZ})$ and $H^1(X, {\bZ})$ (see \cite{Tu} for more details on this rather painful issue). Fix  $T\in H_{0/1}$ such that $\deg (T)=1$, and $K_{0/1}=m_0 T$.   There exists $\chi_0\in H_{0/1}^\sharp$ uniquely determined by the requirements
\[
\chi_0(T)=\rho,\;\;\chi_0\!\mid_{T(H_{0/1})}=1,
\]
where $\rho$ is a primitive $m_0$-th root of $1$. According to  Lemma \ref{lemma: surg} we can think of $\chi_0$ as a character of $G$.\footnote{The Universal Coefficients Theorem and the Poincar\'{e} duality identifies $G^\sharp=H_1(N,\pa N; {\bZ})^\sharp$ with the torsion subgroup of $H_1(N, {\bZ})$.  Via this identification we have $\chi_0={\bf j}\lambda_0\in H_1(N, {\bZ})$.}

\label{ss: 31}

\subsection{Surgery formula for the modified Seiberg-Witten invariant}
We can now state the main surgery formula for the modified Seiberg-Witten invariant. Let $N$, $M_{p/q}$ etc. be as  above. Fix a $spin$ structure $\eps$ on $N$. It induces $spin$-structures $\eps_{p/q}$ on $M_{p,q}$.  For $h\in H_{p,q}$ We will write $\ssw^0_{p/q}(h)$ for $\ssw^0_{M_{p/q}}(h^{-1}\eps_{p/q})$.

\begin{theorem}{\bf (\cite[Marcolli-Wang]{MW}, \cite[Ozsv\'{a}th-Szab\'{o}]{OS})}  For every $p,q$ there exists
\[
f_{p,q,m_0}: U_{m_0}\ra {\bQ}
\]
which  depends only on $p,q,m_0$ but not on $N$ such that for every $g\in G$ we have
\begin{equation*}
\sum_{\pi_{p/q}(h)=g}\ssw^0_{p/q}(h)=p\sum_{\pi_{1/0}(h)=g}\ssw^0_{1/0}(h)+q\sum_{\pi_{0/1}(h)=g}\ssw^0(h)+ f_{p,q,m_0}(\chi_0(g)).
\tag{$\dag_g$}
\label{tag: dag}
\end{equation*}
\label{th: os}
\end{theorem}
 To get more information out of this formula we will  take a partial Fourier transform. Let $\chi\in G^\sharp$.     Using Lemma \ref{lemma: surg} we can identify  $\chi$ with a character of $H_{p/q}$ with the property that $\chi(h)=\chi(h')$ whenever $\pi_{p/q}(h)=\pi_{p/q}$. If we multiply (\ref{tag: dag}) by $\chi(g)$ and then we sum over $g\in G$ we deduce
\[
\widehat{\sw}^0_{p/q}(\chi)=p\widehat{\sw}^0_{1/0}(\chi)+ q\widehat{\sw}^0_{0/1}(\chi) +\sum_{g\in G}f_{p,q,m_0}(\chi_0(g))\chi(g).
\]
To gain further insight we  need to simplify the sum on the right hand side. We have
\[
\sum_{g\in G}f_{p,q,m_0}(\chi_0(g))\chi(g)= \sum_{\rho\in U_{m_0}}\Bigl(\sum_{\chi_0(g)=\rho}f_{p,q,m_0}(\rho)\Bigr)\chi(g)= \sum_{\rho\in U_{m_0}}\Bigl(\sum_{\chi_0(g)=\rho}\chi(g)\Bigr)f_{p,q,m_0}(\rho)
\]
Observe that if $\chi\not\equiv 1$ on $\ker \chi_0$  then
\[
\sum_{\chi_0(g)=\rho}\chi(g)=0,\;\;\forall \rho\in U_{m_0}.
\]
If $\chi\equiv 1$ on $\ker \chi_0$ then   there exists $j\in {\bZ}$ such that $\chi=\chi_0^j$ and
\[
\sum_{\chi_0(g)=\rho}\chi(g)= |\ker \chi_0| \rho^j=\frac{|G|}{m_0}\rho^j.
\]
 Denote by $F_{p,q,m_0}$ the   function
\[
F_{p,q,m_0}: {\bZ}/m_0{\bZ}\ra {\bC},\;\; F_{p,q,m_0}(j\mod {\bZ})=\frac{1}{m_0}\sum_{\rho\in U_{m_0}}f_{p,q,m_0}(\rho)\rho^j.
\]
$F_{p,q,m_0}$ is precisely the Fourier transform of $\frac{1}{m_0}f_{p,q,m_0}$. We  deduce
\begin{equation}
\widehat{\sw}^0_{p/q}(\chi)=p\widehat{\sw}^0_{1/0}(\chi)+ q\widehat{\sw}^0_{0/1}(\chi) +|G|\left\{
\begin{array}{ccc}
F_{p,q,m_0}(j) & {\rm if} & \chi=\chi_0^j\\
& & \\
0 & & {\rm otherwise}
\end{array}
\right. .
\label{eq: surgsw}
\end{equation}

\subsection{Surgery formula for the modified Reidemeister-Turaev torsion}
The modified Reidemeister-Turaev torsion   satisfies a surgery formula very similar in spirit to  (\ref{eq: surgsw}). We  first need to  survey a  few  algebraic facts in the  special setting of the surgery formula. For details and proofs we refer to \cite{N2, Tu2}.

For any finite Abelian group $G$ we set
\[
{\bQ}[G]_0=\Bigl\{ P\in {\bQ}[G];\;\;P(1)=0\Bigr\}.
\]
Consider a rank $1$ Abelian group $A={\bZ}\oplus Tors$, $C$ a finite cyclic group of order $m$, and $\vfi: A\ra C$ a surjective morphism.  Fix a generator  $u\in A$  of  ${\bZ}\subset A$,  and let
\[
{\bZ}[[A]]_+:= {\bZ}[A]+\Theta_A{\bZ}[u,u^{-1}, (1-u)^{-1}].
\]
(We  refer to \cite{Tu2,Tu} for an invariant definition of ${\bZ}[[A]]_+$, which does not rely on the non-canonical decomposition $A=  \mbox{free part}\oplus \mbox{torsion}$.) The morphism $\vfi$ induces a morphism
\[
\vfi_\sharp:{\bZ}[[A]]_+\ra {\bZ}[C]_0.
\]
Its definition is best expressed in terms of Fourier transforms.   Think of an element $f\in {\bZ}[[A]]_+$ as a function $f:A\ra {\bZ}$.  As such, it has a Fourier transform
\[
\hat{f}: A^\sharp \ra {\bC},\;\;\hat{f}(\chi)= \sum_{a\in h} f(a)\chi(h).
\]
The above infinite sum  may not be convergent for all $\chi$, but the  Fourier  transform still makes sense as a distribution on the compact Lie group $A^\sharp$.  There are finitely many characters $\chi$ for which the series is not convergent, and they are  characterized by the condition $\chi(u)=1$.  We will call such characters {\em singular} and the other regular. We denote by $A^\sharp_{reg}$ the set of regular characters, and set $A^\sharp_{sing}:=A^\sharp\setminus A^\sharp_{reg}$. The morphism $\vfi$  induces by duality an inclusion
\[
\vfi^\sharp: C^\sharp\hookrightarrow A^\sharp.
\]
The Fourier transform of $\vfi_\sharp f$ is a function $\widehat{\vfi_\sharp f}: C^\sharp\ra{\bC}$. Let $\chi\in C^\sharp$. In \cite{N2} we have shown how to compute the value $\widehat{(\vfi_\sharp f)}(\chi)\in {\bC}$.  More precisely we have shown the following.

\medskip

\noindent $\bullet$   If $f\in {\bZ}[A]$, so that  $f$ has  finite support as a function $A\ra {\bZ}$, we have
\begin{equation}
\widehat{(\vfi_\sharp f)}(\chi)=\left\{
\begin{array}{ccc}
\hat{f}(\vfi^\sharp\chi) & {\rm if} & \chi\neq 1\\
0 & {\rm if} & \chi =1
\end{array}
\right. .
\label{eq: fourier1}
\end{equation}
$\bullet$ If $f=\Theta_A(1-u)^{-1}$ then
\begin{equation}
\widehat{(\vfi_\sharp f)}(\chi)=\left\{
\begin{array}{ccc}
\widehat{\Theta}_A(\vfi^\sharp\chi)\bigl(1-\vfi^\sharp\chi(T)\bigr)^{-1}& {\rm if} & \vfi^\sharp\chi\in  A^\sharp_{reg}\\
& &\\
0 & {\rm if} & \vfi^\sharp\chi \in   A^\sharp_{sing}
\end{array}
\right. .
\label{eq: fourier2}
\end{equation}
Suppose  now that $\chi: A\ra U_m$ is a surjective character, and $f\in {\bZ}[[A]]_+$.   The identity function $\iota_m: U_m\ra U_m$ is a character of $U_m$ and $\chi^\sharp(\iota_m)=\chi$. We get an element $\chi_\sharp f\in {\bZ}[U_m]_0$. If $f\in {\bZ}[A]$ then
\begin{equation}
(\chi_\sharp f)(\iota_m)=\left\{
\begin{array}{ccc}
\hat{f}(\chi) &{\rm  if}  & m>1\\
0 & {\rm if}& m=1
\end{array}
\right..
\label{eq: fourier3}
\end{equation}
If $f= \frac{\Theta_A}{1-u}$ then
\begin{equation}
\widehat{(\vfi_\sharp f)}(\iota_m)=\left\{
\begin{array}{ccc}
\widehat{\Theta}_A(\chi)\Bigl(1-\chi(T)\Bigr)^{-1}& {\rm if} & \chi(T)\neq 1\\
& &\\
0 & {\rm if} & \chi(T)=1
\end{array}
\right. .
\label{eq: fourier4}
\end{equation}
If $\vfi:A\ra B$ is a surjective morphism of finite Abelian groups then we get  morphisms
\[
\vfi_\sharp: {\bQ}[A]_0\ra {\bQ}[B]_0,\;\;\vfi^\sharp: B^\sharp\ra A^\sharp.
\]
Then for $f\in {\bQ}[A]_0$ and $\chi\in B^\sharp$ we have
\[
\widehat{(\vfi_\sharp f)}(\chi)=\left\{
\begin{array}{ccc}
\hat{f}(\vfi^\sharp\chi) & {\rm if} & \chi\neq 1\\
0 & {\rm if} & \chi =1
\end{array}
\right. .
\]
We can now return to topology. We will continue to use the notations in the previous section. Applying \cite[Lemma 6.2]{Tu} iteratively we deduce the following result.

\begin{theorem} Suppose $\chi $ is a nontrivial character of     $G=H_1(N,\pa N; {\bZ})$, so that $\chi(G)=U_m$, for some $m>1$. Then
\[
\chi_\sharp\t_{p/q}=p\chi_\sharp\t_{1/0}+q\chi_\sharp\t_{0/1}.
\]
\label{th: torsurg}
\end{theorem}
Above and in the sequel  we  use the convention $Object_{p/q}:=Object(M_{p/q})$. To proceed further we  take the Fourier transform of the above formula and we get
\[
\widehat{(\chi_\sharp\t_{p/q})}(\iota_m)=\widehat{(p\chi_\sharp\t_{1/0})}(\iota_m)+\widehat{(q\chi_\sharp\t_{0/1})}(\iota_m),
\]
where  $m={\rm ord}\,(\chi)$. Recall that
\[
\t_{p/q}=\t_{p,q}^0-\frac{1}{2}CW_{p/q}\Theta_{p/q},\;\; \t_{1/0}=\t^0_{1/0}-\frac{1}{2}CW_{1/0}\Theta_{1/0},
\]
\[
\t_{0/1}= \t_{0/1}^0+\frac{\Theta_{0/1}T}{(1-T)^2}.
\]
Using (\ref{eq: fourier3}) and (\ref{eq: fourier4}) and the identities $\Theta_{p/q}(\chi)=0=\Theta_{1/0}(\chi)$ for $\chi\neq 1$ we deduce we deduce
\[
\widehat{\t}_{p/q}^0(\chi) =p\widehat{\t}^0_{1/0}(\chi)+ q\widehat{\t}_{0/1}^0(\chi)+\left\{
\begin{array}{ccc}
\widehat{\Theta}_{0/1}(\chi)\chi(T)\Bigl(1-\chi(T)\Bigr)^{-2}& {\rm if} & \chi(T)\neq 1\\
& &\\
0 & {\rm if} & \chi(T)=1
\end{array}
\right. .
\]
The last term is nontrivial if and only if $\chi(T)\neq 1$ and $\chi\mid_{T(H_{0/1}}=1$. This is possible if and only if $\chi=\chi_0^j$, for some $j=1,2,\cdots, m_0-1$. Additionally $\Theta_{0/1}(\chi_0^j)=|T(H_{0/1}|=|G|/m_0$. We conclude that if $\chi$ is a nontrivial character of $G$ we have
\[
\widehat{\t}_{p/q}^0(\chi) =p\widehat{\t}^0_{1/0}(\chi)+ q\widehat{\t}_{0/1}^0(\chi)+\frac{|G|}{m_0}\left\{
\begin{array}{ccc}
\frac{\chi_0^j}{(1-\chi_0^j)^2} & {\rm if} & \chi=\chi_0^j,\;\;j=1,\cdots, m_0-1\\
& &\\
0  & {\rm if} & \chi\neq \chi_0^k,\;\;k=1,\cdots ,m_0
\end{array}
\right.
\]
 We need to relate
 \[
\widehat{ \t}^0_{p/q}(1)=\frac{|H_{p/q}|}{2}CW_{p/q}=\frac{pm_0|G|}{2}CW_{p/q},
 \]
 \[
\widehat{\t}^0_{1/0}(1)=\frac{|H_{1/0}|}{2}CW_{1/0}= \frac{m_0|G|}{2}CW_{1/0},
 \]
 and
 \[
\widehat{\t}^0_{0/1}(1)=\frac{1}{2}\Delta_{M_{0/1}}''(1).
 \]
 This follows from  the surgery formula for the Casson-Walker invariant \cite[\S 4.6]{Lescop}, \cite[Chap.4]{Walker}.  More precisely,  the arguments in \cite[p.38-39]{OS} yield
 \[
 CW_{p/q}= pCW_{1/0} +\frac{q}{2}\Delta''_{M_{0/1}}(1) + |G|\Biggl(\frac{q(m_0^2-1)}{12m_0}-\frac{pm_0s(q,p)}{2}\Biggr),
 \]
where $s(q,p)$ denotes the Dedekind sum of the pair of  co-prime integers $q,p$.  Putting all of the above together we deduce that  for every pair of coprime integers $(p,q)$ and every positive integer $m_0$ there exists a function $G_{p,q,m_0}: {\bZ}/m_0{\bZ}\ra {\bC}$ such that
\begin{equation}
\widehat{\t}_{p/q}^0(\chi) =p\widehat{\t}^0_{1/0}(\chi)+ q\widehat{\t}_{0/1}^0(\chi)+|G|\left\{
\begin{array}{ccc}
G_{p,q,m_0}(j) & {\rm if} & \chi=\chi_0^j\\
& & \\
0 & & {\rm otherwise}
\end{array}
\right. .
\label{eq: surgtors}
\end{equation}
The similarity with (\ref{eq: surgsw}) is striking.    The results in \cite{MW, OS} show that
\[
F_{p,q,m_0}(1)=G_{p,q,m_0}(1),\;\;\;\forall p,q, m_0.
\]
In particular
\[
F_{p,q,1}=G_{p,q,1},\;\;\forall p,q.
\]
Let us briefly comment on the ``flavor'' of the surgery
formul{\ae} (\ref{eq: surgsw}) and (\ref{eq: surgtors}).      Note
first that  the   first homology group of a rational homology
$3$-sphere can be naturally identified with its dual using the
linking form.   We can think of the invariants $\t^0_M$ and
$\sw_M^0$  as functions on $H^1(M, {\bZ})$, as well as functions on the dual.

Suppose we perform Dehn surgery on a knot $K\hra M$ to obtain a new rational homology sphere $M(K)$. The surgery formula     essentially states that  if we know the values of these invariants on  homology classes $c\in H_1(M, {\bZ})$ which do not link with $K$ then we can also compute the values of these invariants on   homology classes $c\in H_1(M(K), {\bZ})$ which do not link with  $K\hra M(K)$.

More rigorously, consider a pair $M_0$, $M_1$  related by a  Dehn surgery on a knot $K$. Denote by $N$ the common knot complement, and set
\[
H:= H_1(N, {\bZ}), \;\;G:= H_1(N,\pa N; {\bZ}),\;\;H_i:= H_1(M_i,{\bZ}),\;\;i=0,1.
\]
We have a diagram of surjective morphisms
\[
\begin{diagram}
\node{}\node{H}\arrow{sw,t,A}{\pi_0}\arrow{se,t,A}{\pi_1}\node{}\\
\node{H_0}\node{}\node{H_1}
\end{diagram}
\]
Dualizing we get the  diagram
\[
\begin{diagram}
\node{}\node{H^\sharp}\node{}\\
\node{H^\sharp_0}\arrow{ne,t,J}{\pi^\sharp_0}\node{}\node{H^\sharp_1}\arrow{nw,t,L}{\pi^\sharp_1}
\end{diagram}
\]
The group $G$ can be identified with a subgroup of  $H^\sharp$. The knot $K$ defines two subgroups
\[
K_i^\perp :=\Bigl\{ c\in H_i^\sharp;\;\; c(K_i)=1\Bigr\}=\Bigl\{ c\in H_i^\sharp; \;\;{\bf lk}_{M_i}(c, K_i)=0\Bigr\},\;\;i=0,1,
\]
and we have isomorphisms $\pi_i^\sharp: K_i^\perp\ra G$. We  can think of $G$  as the graph of a correspondence $T_K: H_0^\sharp\ra H_1^\sharp$ induced by the Dehn surgery.  The domain of this correspondence is $K_0^\perp$,  the range is $K_1^\perp$, and viewed as a  correspondence $T_K\subset K_0^\perp\times K_1^\perp$ it is a group isomorphism.    We will refer to a such a correspondence  as a {\em partial isomorphism} (p.i.) of groups.

For a surgery  along a  knot $K\hra M$,  whose meridian satisfies $\lambda_0\cdot \mu=1$, and attaching curve $c=p\mu+q\lambda_0$, we will denote by $\xi:=\xi_{K,c}$ the induced  p.i.   We denote by  $G_K$ the group $H_1(M\setminus K,\pa(M\setminus K); {\bZ})$, by $m_0(K,c)$ respectively $p=p(K,c)$ the  divisibility,  and respectively  multiplicity of the surgery. Finally, set $D_M:=\sw^0_M-\t^0_M$.  Since $\t^0_M=\sw^0_M$ if $b_1(M)>0$ we can now rephrase the surgery formul{\ae} (\ref{eq: surgsw}) and (\ref{eq: surgtors})
\[
\widehat{D}_{M_{k,c}}( \xi_{K,c}\chi)= p \widehat{D}_{M} +|G_K|\k_{p,q,m_0},\;\;\forall \chi\in K^\perp= {\rm Dom}\,(\xi_{K,c}),
\]
where the correction term is a function  on $G_K$, which is nontrivial only on  the cyclic group of order $m_0$ generated by ${\bf j}\lambda_0\in T(H_1(N\setminus K, {\bZ}))=G_K^\sharp$, and depends only on the arithmetic of the surgery, $p,q,m_0$. Moreover $\k_{p,q,1}=0$.

More generally  consider a  $3$-manifold $N$ such that $\pa N\cong
T^2$, fix a longitude $\lambda_0$, and two primitive classes $c_0,
c_1$  represented by two simple closed curves.       By Dehn
surgery with attaching curves $c_0,c_1$ we get two manifolds
$M_{c_0}$, $M_{c_1}$, with first homology groups $H_{c_0}$,
$H_{c_1}$, and  distinguished classes $K_{c_i}\in H_{c_i}$,
$i=0,1$, defined by the core of the attached solid torus.  Set
$G:=H_1(N,\pa N; {\bZ})$, and  denote by $\xi_{c_1,c_0}$  the
isomorphism $K_{c_0}^\perp\ra K_{c_1}^\perp$ induced by the
surgery cobordism.  We denote by $[c_0,c_1]$  the orbit of  $(c_0,c_1)$ relative to  the action of $SL_2({\bZ})$ on the space of pairs of primitive classes $c_0,c_1\in H_1(\pa N, {\bZ})$. Then we have
\begin{equation}
(\lambda_0\cdot c_0)\widehat{D}_{M_{c_1}}(\xi_{c_0,c_1}\chi)=(\lambda_0\cdot c_1)\widehat{D}_{M_{c_0}}(\chi) + |G|\k_{[c_0,c_1],m_0}(\chi),\;\;\forall \chi\in K_{c_0}^\perp.
\label{eq: surg}
\end{equation}
The {\bf arithmetic type} $\alpha$ of a surgery is  the  pair $([c_1,c_2], m_0)$.   We  denote by $\A$ the   set of all arithmetic types  for which the correction term $\k$ is trivial.  We know that
\[
([c_1,c_2], 1)\in \A,\;\;\forall c_1,c_2.
\]
We will call the surgeries of arithmetic type $\alpha \in \A$ as {\bf admissible}.

\begin{remark} As explained in \cite[Remark B.6]{N2}, the orbit  $[c_0,c_0]$ is  completely characterized by  the  extension
\[
0\ra {\bZ}\lan c_0\ran\oplus {\bZ}\lan c_1\ran\hra H^1(\pa N, {\bZ})\ra H^1(\pa N, {\bZ})/\Bigl( {\bZ}\lan c_0\ran\oplus {\bZ}\lan c_1\ran\Bigr)\ra 0.
\]
More precisely, the quotient group  $ H^1(\pa N, {\bZ})/( {\bZ}\lan c_0\ran\oplus {\bZ}\lan c_1\ran)$ is isomorphic  to  the cyclic groups of order $|c_0\cdot c_1|$, and the extension is characterized by a  character of this group. Thus, the orbit $[c_0,c_1]$  is described by the integer $c_0\cdot c_1$,   and a character of   ${\bZ}_{|c_0\cdot c_1|}$.
\end{remark}

\section{Seiberg-Witten $\Longleftrightarrow$ Casson-Walker+ Reidemeister  torsion}
\setcounter{equation}{0}
\subsection{Topological preliminaries}
\label{ss: 41}
 Denote by $\X$ the family of all closed, compact oriented $3$-manifolds $M$ such that
\[
\sw_M^0=\t_M^0.
\]
We want to prove that $\X$ consists of all $3$-manifolds.

We already know that $M\in X$ if  $b_1(M)>0$, or  $M$ is an {\em integral} homology sphere, or if $M$ is a lens space.  Also, we have
\[
M_1,M_2\in \X \Longrightarrow M_1\# M_2\in \X.
\]
\begin{definition} A {\bf deflating primitive surgery}  is a Dehn surgery on a knot $K$ in a rational homology sphere  $M$ with the following  properties.

\noindent (a) The longitude  $\lambda \in H_1(\pa M\setminus K, {\bZ})$ is a primitive class.

\noindent (b) The attaching curve $c$ of the surgery satisfies $c\cdot \lambda=\pm 1$.

An {\bf excellent} surgery is a deflating primitive surgery  which does not change the order of the first homology group.  Two rational homology $3$-spheres  will be called {\bf e-related} if one can be obtained from the other by a sequence of excellent surgeries.
\end{definition}
The attribute deflating is justified by the inequality
\[
|H_1(M', {\bZ})|\leq |H_1(M, {\bZ})|
\]
when $M'$ is obtained from $M$  by a deflating primitive surgery. The surgery is excellent  iff we have equality.    Note that if $M_0$ and $M_1$ are e-related  then $H_1(M_0, {\bZ})\cong H_1(M_1, {\bZ})$ and  they have isomorphic linking forms. The following result is immediate.

\begin{lemma}
Suppose  $M\in \X$, $b_1(M)=0$.   If $M'$ is obtained from $M$ by a deflating primitive surgery then $M'\in \X$. In particular, if $M\in \X$ and $M'$ is e-related to $M$ then $M'\in \X$.
\label{lemma:1}
\end{lemma}
\noindent{\bf Proof}\hspace{.3cm}   Indeed, we have $G:=H_1(M',{\bZ})\cong H_1(M\setminus K, \pa (M\setminus K); {\bZ})$  and $F_{p,q,1}=G_{p,q,1}$. The surgery formul{\ae}   establish the equality of $\t^0_{M'}$ and $\sw^0_{M'}$ as functions on $G^\sharp$, and $G^\sharp$ turns out to be their  maximal domain. $\Box$

\medskip

\begin{corollary} $\X$ contains lens spaces, integral surgery spheres, and is closed under connected sums and deflating  primitive surgeries.
\label{cor: ratl}
\end{corollary}

Before we proceed further  we want to  briefly recall some basic topological facts. For more details we refer to \cite{GS, Sav}. Any   rational homology sphere can be  obtained from $S^3$ by performing Dehn surgery on a link $L= K_1\cup \cdots \cup K_n$  with  surgery coefficients  $p_1/q_1,\cdots ,p_n/q_n$.  Set $\vec{p}=(p_1,\cdots, p_n)$, $\vec{q}=(q_1,\cdots ,q_n)$. We denote by $M=M(L,\vec{p},\vec{q})$  the three manifold obtained by this surgery.  We say that a surgery diagram belongs to $\X$ if the corresponding  $3$-manifold belongs to $\X$.

All the homological data  of $M(L,\vec{p},\vec{q})$ is contained in the $n\times n$ linking matrix $A=A(L; \vec{p},\vec{q})$ defined by
\[
a_{ij}=\left\{
\begin{array}{ccc}
p_i/q_i &{\rm if}  & i=j\\
\ell_{ij} & {\rm if} & i\neq j
\end{array}
\right. ,
\]
where $\ell_{ij}={\bf Lk}(K_i,K_j)$.  The surgery diagram is called integral if   the linking matrix is integral.

Denote by $\mu_i$ the meridian of $K_i$, set $\Omega:=A^{-1}$, and denote by $Q$ the diagonal $n\times n$ matrix $Q:= {\rm diag}\,(q_1,\cdots,q_n)$. The manifold  $M$ is a rational  homology sphere if and only if the linking matrix $A$ is nonsingular.  The first homology group admits the presentation
\[
0\ra {\bZ}^n \stackrel{QA}{\Lra} {\bZ}^n \ra H_1(M, {\bZ})\ra 0
\]
so that its  order is $\det (QA)$. Moreover the  images of the knots $K_i$  generate $H_1(M, {\bZ})$  and  we have
\[
{\bf lk}_M(K_j,K_i)=-\Omega_{ji}\mod {\bZ}.
\]
We  have a natural isomorphism $H_1(S^3\setminus L ;{\bZ})\ra {\bZ}^n$ defined by
\[
c\mapsto \Bigl(\,{\bf Lk}(c,K_1),\cdots, {\bf Lk}(c,K_n)\,\Bigr)
\]
for any closed curve disjoint from $L$.  More geometrically
\[
c=\sum_{i=1}^m{\bf Lk}(K_i,c)\mu_i.
\]
Such a closed curve  defines a homology class $[c]$. We have $[\mu_i]=-q_i[K_i]$.

Suppose $[c]$ is a homology class  in  $M$ of order  $m$. (We set $m=1$ if $[c]=0$.) A  surgery on a knot  representing $[c]$ has divisibility $m_0$ determined   by
\[
m_0:= (k,m),\;\;{\bf lk}_M([c],[c])=\frac{k}{m}\mod {\bZ}.
\]
A class $c\in H_1(M, {\bZ})$  is called {\em primitive} if it has divisibility one.   Note that if $K$ is a node in $M$ representing a primitive class of order $m$,   and $M'$ is a manifold obtained from $M$ by deflating surgery  then
\[
|H_1(M', {\bZ})|=\frac{1}{m}|H_1(M, {\bZ})|.
\]
The above observations show  that the  excellent  surgeries  are precisely the $1/q$-surgeries on a  homologically trivial knots.

The  {\bf pruning}  of a surgery diagram is the operation  of removing the components with surgery coefficients $\pm 1$ which are algebraically split from the rest of the diagram. The  pruning   is equivalent to performing a sequence of excellent    surgeries. We    say that two  surgery diagrams are {\bf p-related} if  one can go form one to another by a sequence of Kirby moves and prunings.

\begin{corollary} If $\mathcal{D}$ is a surgery diagram p-related to a diagram in $\X$ then  $\mathcal{D}$ is also in $\X$.
\end{corollary}

For every $\vec{a}\in {\bZ}^n$ we denote by $[\vec{a}]$ the rational number
\[
[\vec{a}]=a_1-\frac{1}{a_2-\frac{1}{a_3-\ddots}}.
\]
\begin{lemma}{\bf \cite[N. Saveliev]{Sav1}} Any homology lens space is e-related to a lens space.
\label{lemma: sav}
\end{lemma}

\begin{figure}[ht]
\begin{center}
\epsfig{figure=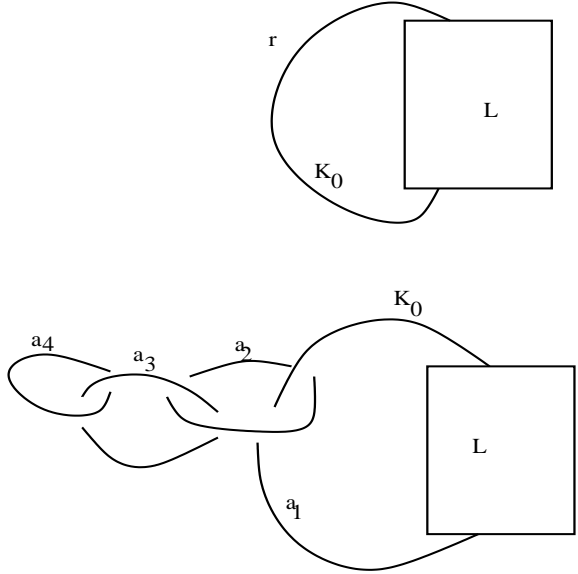, height=3in, width=3in}
\end{center}
\caption{\sl {Slam-dunking $K_0$.}}
\label{fig: 1}
\end{figure}

\noindent {\bf Proof} \hspace{.3cm}  Any homology lens space $M$ is  obtained by Dehn surgery  on a knot $K_0$  in an integral homology sphere $M'$, \cite{BL}. Denote by $r\in {\bQ}$ the surgery coefficient of $K_0$.   We can represent the homology sphere $M'$  as surgery on an algebraically split   link $L=K_1\cup \cdots \cup K_n$ in $S^3$ with surgery coefficients  $\ve_j=\pm 1$. We can think of $M$ as obtained from $S^3$ by surgery on the link $L_0=K_0\cup L$.     Suppose
\[
r=[\vec{a}],\;\;\vec{a}=(a_1,a_2,\cdots, a_m).
\]
Performing a sequence of slam-dunks as in \cite[Sec. 5.3]{GS}  we can replace   $L_0$ with the link $L\cup K$  as in Figure \ref{fig: 1}.
We have thus succeeded in presenting   $M$ as  an integral surgery on a link in $S^3$ with linking matrix
\[
\left[\begin{array}{ccccccccccc}
\pm 1 & 0 & 0 & \cdots & 0                     & & \ell_1 &0 &\cdots &\cdots &0\\
0 & \pm 1 & 0 & \cdots & 0                     &  & \ell_2 & 0 &0 &0 &0 \\
\vdots & \vdots & \vdots & \vdots & \vdots     & & \vdots & \vdots & \vdots & \vdots & \vdots \\
0   & 0 & \cdots & \cdots    & \pm 1           &     &\ell_n & 0 &\cdots &\cdots &0\\
&&&&&&&&&& \\
\ell_1 & \ell_2 &\cdots & \cdots & \ell_n      & & a_1   & 1 & 0 & 0 & \cdots\\
0 & 0 & \cdots & \cdots           & 0          & & 1 & a_2 & 1 & 0 &\cdots\\
\vdots & \vdots & \vdots & \vdots & \vdots     & &  0 & 1 & a_3 & 1 &   \cdots\\
\vdots & \vdots & \vdots & \vdots & \vdots     & &\vdots & \vdots & \vdots & \vdots & \vdots \\
0 &\cdots &  \cdots & \cdots & 0                & &  0 & \cdots & \cdots & 1& a_m
\end{array}
\right]
\]
The first part of this matrix is described by the link $L$, and $\ell_j:= {\bf Lk}(K_0, K_j)$. By sliding $K_0$ over the components  of $L$ we  can   kill the off-diagonal terms $\ell_i$. This changes the topological type of $K_0$. It becomes a knot $K_0'$,  and the surgery coefficient $a_1$ changes to some integer $a_1'$.   Reverse the slam-dunks. We get a new link  algebraically split link $L_2= K_0'\cup L$ where the surgery coefficient of $K_0'$ is  $r'=[a_1',a_2,\cdots, a_m]$, and the surgery  coefficient of $K_j$ is $\pm 1$. By inserting  $\infty$-unknots and performing a sequence of Rolfsen twists we can replace $K_0'$ with an unknot $K_0''$.  We  can thus describe $M$ as surgery on the algebraically split link
\[
L_2=K_0''\cup K_1\cdots \cup K_n
\]
with surgery coefficients $\ve_0=r'$, $\ve_1=\pm 1,\cdots,\ve_n=\pm 1$.  The $r'$ surgery on $K_0''$ is a lens space while the surgeries on $K_j$ are excellent surgeries. This shows $M$ is e-related to a lens space. $\Box$

\subsection{Proof of the main result}
 We will present a proof  by induction over the ``complexity'' of a  rational homology sphere. To define the notion of complexity we need to   present a few algebraic facts about the linking forms of such manifolds. We follow the notations in \cite{KK}.

 For each prime $p>1$, each $q\in {\bZ}$, and each $k\geq 1$ such that $(p,q)=1$ denote by $A_p^k(q)$ the  linking form n the cyclic group ${\bZ}/p^k$  defined by
 \[
 {\bf g}\cdot {\bf g}=\frac{q}{p^k}
 \]
 where $g$  denotes the natural generator of this group.     Also denote by $E_0^k$ , $k\geq 1$ and $E_1^k$, $k\geq 2$, the linking forms on ${\bZ}/2^k\oplus {\bZ}/2^k$ defined by the matrices
 \[
 E_0^k=\left[
 \begin{array}{cc}
 0 & 2^{-k}\\
 2^{-k} & 0
 \end{array}
 \right],\;\;E_1^k=\left[
 \begin{array}{cc}
 2^{1-k} & 2^{-k}\\
 2^{-k} & 2^{1-k}
 \end{array}
 \right].
 \]
When referring to $A_p^k(q)$, $E_0^k$ , $E_1^k$ we  mean the corresponding groups equipped with  these linking forms. Define the complexity of $A_p^n$ to be $\kappa(A^n_p)=p^{k+1}$. Define  the complexity of $E_i^n$, $i=0,1$ to be $\kappa(E^n_i)=2^{2n+2}$.

A classical result of C.T.C. Wall \cite{Wall}  shows that every
linking form $(G,{\gq})$  decomposes non-uniquely into an
orthogonal sum of $A$'s and $E$'s.      If ${\gq}$ is a  linking
form on a $p$-group, then we define   its complexity to  be  the
product of the complexities of its  elementary  constituents $A$ and/or $E$  in some orthogonal decomposition. The  results in \cite{KK}  show that this number is
independent of   the chosen  orthogonal decomposition of ${\gq}$ in elementary parts  $A$ and $E$. We denote by $\kappa({\gq})$ the complexity  of
${\gq}$. For every ${\bQ}HS$ $M$ we denote  by $\nu_M$ the order
of $H_1(M, {\bZ})$  , by ${\gq}_M$ the linking form of $M$, and by
$\kappa_M$ the complexity of ${\gq}_M$.  We have  the following
elementary result whose proof is left to the reader.

\begin{lemma} If $M_1$ and $M_2$  are two rational homology spheres such that $\nu_{M_1}|\nu_{M_2}$ and $\nu_{M_1}<\nu_{M_2}$ then $\kappa_{M_1}<\kappa_{M_2}$.
\label{lemma: c1}
\end{lemma}
We would like to present a few    methods of reducing the complexity of a manifold.  The primitive deflating surgeries provide one first example.

\begin{definition}  Let $K$ be a knot  in a rational homology sphere $M$ supporting a {\bf nontrivial} homology class. The knot $K$ is called {\bf good} if ${\gq}_M(K,K)\neq 0$. Otherwise, it is called {\bf bad}.
\end{definition}
Suppose $K$ is a good knot in a rational homology sphere  $M$. If $r$ is order of $K$ then
\[
{\gq}(K,K)= \frac{m}{r},\;\;0<m<r,
\]
and the divisibility of any surgery on this knot is $m_0(K):=(m,r)$.      Consider any surgery with attaching curve $c$ satisfying $|c\cdot \lambda|= m_0$. This is a surgery of divisibility $m_0$ and of type $(p,q)=(m_0,\ast)$. We obtain  a new rational homology sphere $M'$ such that $\nu_{M}= \frac{r}{m_0}\nu_{M'}$. Lemma \ref{lemma: c1} shows that  the complexity of $M'$ is smaller than the complexity of $M$.  We have thus proved the following result.

\begin{corollary}         The complexity of a  rational homology sphere can be reduced by performing surgeries on  good knots. Moreover if  the original manifold has no $2$-torsion we can arrange so that the resulting manifold also has no $2$-torsion.
\label{cor: c1}
\end{corollary}

Certain surgeries on certain bad knots also do reduce the  complexity.    We have  the following technical result whose proof  is deferred to  an Appendix.

\begin{lemma}  Suppose $M$ is  a  rational homology sphere  with linking form  $A^s_{p}(q_1)\oplus A^r_{p}(q_2)$ $s\geq r$, and $K$ is a bad knot in $M$ of the form
\[
K= c_1\oplus c_2
\]
where  the homology class $c_2$ generates $A_{p^r}(q_2)$. Then one can perform a surgery  on $K$ such that the resulting  manifold is a homology lens space of the same order as $M$.
\label{lemma: c2}
\end{lemma}
We will   call the surgery in this lemma $A_p$-surgery. A knot  with  the properties  in the lemma will be called a {\bf mildly bad} knot. Set
\[
\q:=\Bigl\{ {\gq};\;\;{\gq}_M\cong {\gq}\Longrightarrow M\in \X\Bigr\}.
\]
We already know that all the linking forms $A_p^k(q)$ belong to  $\q$.

We need to talk a little bit about admissible surgeries, i.e. surgeries for which  the correction term in the surgery formula (\ref{eq: surg}) is trivial.   Observe that if two rational homology spheres in $\X$ are related by a Dehn surgery then this  surgery is admissible.
\begin{corollary}  The surgeries $A_p$ surgeries described in Lemma \ref{lemma: c2} are admissible.
\label{cor: c2}
\end{corollary}

\noindent {\bf Proof}\hspace{.3cm} Consider a direct sum of two
lens spaces with the above linking forms.  This is a manifold in $\X$. The result of this
surgery produces a rational homology space  which is also a manifold in
$\X$ so the surgery is admissible. $\Box$

\bigskip

 We also want to mention the following topological result. For a proof we refer to \cite{N2}.
\begin{lemma}  Suppose  $M_1$, $M_2$  are two rational homology spheres and $\phi$ is an isomorphism
\[
\phi:(H_1(M_1, {\bZ}), {\gq}_{M_1})\ra (H_1(M_2, {\bZ}), {\gq}_{M_2}).
\]
Suppose $K_i$ is a knot in $M_i$, $i=1,2$ such that $\phi([K_1]) =[K_2]$.  If  $M_i'$ $i=1,2$, are obtained  perform surgeries {\bf of the same  arithmetic type} $\alpha$ on $K_1$ and $K_2$, then there exists an isomorphism
\[
\phi_\alpha:(H_1(M_1', {\bZ}), {\gq}_{M_1'})\ra (H_1(M_2', {\bZ}), {\gq}_{M_2'}).
\]
\label{lemma: c3}
\end{lemma}
The main trick used  in the proof is the following immediate consequence of the surgery formula (\ref{eq: surg}).

\begin{lemma}  Suppose $M$  is a rational homology sphere, and $\chi$ a character  of $H=H_1(M, {\bZ})$.   We identify $H$ with its dual using the linking form.    Suppose that   there exists an admissible surgery  on a knot $K_\chi$ such that
\[
{\gq}_M(\chi, K)=0\Longleftrightarrow\chi\in K^\perp
\]
and the result of the surgery  is a manifold in $\X$. Then $\widehat{D}_M(\chi)=0$. In particular, if for every $\chi$ there exists a knot with the above properties then $\widehat{D}_M\equiv 0$.
\label{lemma: c4}
\end{lemma}

The proof of Theorem \ref{th: main} will be carried out in several steps.

\bigskip

\noindent {\bf Step 1}  Fix a prime number $p>2$, and denote by $\R_p$ the family of rational homology spheres such that$\nu_M=p^r$, $r>0$. We will show  that $\R_p\subset\X$.  The proof will be an induction on the complexity.  For $\kappa\ge0$ denote by $\R^\kappa
$ the manifolds in $\R_p$ of complexity $\leq \kappa$.

The  rational homology lens spaces of order $p$  have  minimal nonzero complexity $p+1$, and belong to $\X$ so that
\[
\R_p^{p+1}\subset \X.
\]
Observe that if $M\in \R_p$ then we have a decomposition
\[
{\gq}_M=\bigoplus_{j=1}^n A_p^{s_j}(q_j),\;\; 0<s_1\leq s_2\leq \cdots \leq s_k,
\]
and
\[
 \kappa_M =  p^{s_1+\cdots s_k+k}.
\]
The integer $k$ is called the {\em rank}, and we denote it  by $\rho_M$.   Define the  {\bf the standard model} of $M$ to be the connected sum  of lens spaces with the same  linking form as $M$.   We denote the standard model   by $\tilde{M}$. Note that for every  $M\in \R_p$ we have $\tilde{M}\in \X$.

Suppose   $\R_p^{\kappa}\subset \X$. We want to prove  that $\R_p^{\kappa+1}\subset \X$.  Let $M\in \R_p^{\kappa +1}$. Set $H:=H_1(M, {\bZ})$, and  fix a nontrivial character $\chi$ of $H$.     We distinguish two cases.

\medskip

\noindent {\bf Case 1} There exists  a good knot $K\in \chi^\perp$.     Then there exists a good knot $\tilde{K}$ in the model   $\tilde{M}$. We can perform   a complexity reducing surgery on  $\tilde{K}$ to obtain a manifold of  smaller complexity  which by induction we know is in $\X$. This show that the  arithmetic type of this surgery is admissible.  We perform this admissible surgery on the   knot $K$ on $M$ and we obtain  a manifold of smaller complexity.   Lemma \ref{lemma: c4} then implies $\widehat{D}_M(\chi)=0$.

\medskip

\noindent {\bf Case 2} $\chi^\perp$ consists only of bad knots.   If the rank of $H$ is $1$ then $M$ is a rational homology space so it is in $\X$.  Suppose the rank is $>1$.      $(H,{\gq}_M)$ decomposes into a nontrivial sum of cyclic $p$ groups
\[
H= {\bZ}/p^{s^1}\oplus \cdots \oplus {\bZ}/p^{s^k} ,\;\; 0<s_1\leq \cdots \leq s_k,\;\; k >1.
\]
We get a corresponding decomposition  $\chi=\chi_1\oplus \cdots \oplus \chi_k$. Observe that all components must be nonzero. Indeed, if $\chi_j=0$ then  the generator of  the $j$-th component belongs to  $\chi^\perp$ and is  a  good knot.  Thus $\chi_1,\chi_2\neq 0$.    It is easy to see that $\chi^\perp \cap {\bZ}/p^{s^1}\oplus {\bZ}/p^{s_2}\neq 0$.  Pick a mildly bad knot $K$ in this  group. Thus all but the first two components of $K$ are zero, and one of the components generates the corresponding summand. Perform an  $A_p$ surgery on this  knot.  This reduced the complexity of $M$. By induction, the resulting manifold is in $\X$ that $\hat{D}_M(\chi)=0$. Thus $\R_p\subset \X$

\bigskip

\noindent {\bf Step 2} If $\nu_M$ is odd then $M\in \X$.   For each vector $\vec{p}=(p_1,\cdots, p_n)$ whose components consist of  pairwise distinct of odd primes. Denote by $\R_{\vec{p}}$  the family  of  rational homology spheres $M$ such that the prime divisors of $\nu_M$ are amongst the primes $p_j$.  Again we perform induction on   complexity. The  first homology group $H$ of $M\in \R_{\vec{p}}$ decomposes  as an orthogonal direct sum of $p$-groups
\[
H=\bigoplus_{j=1}^nG_{p_j},\;\;|G_{p_j}|=p_j^{s_j}.
\]
Each component $G_{p_j}$ decomposes  as an orthogonal sums of $A^\ast_{p_j}(\ast)$. Denote by $r_j$ the number of such components.   We have
\[
\kappa(M)= |H|+ r_1+\cdots + r_n.
\]
Suppose $\chi$ is a nontrivial character of $H$.   We distinguish again two cases.

\medskip

\noindent {\bf Case 1} $\chi^\perp$ contains   good knots.  In this case  we  perform  a surgery  as in  Lemma \ref{cor: c1} which produces a manifold of smaller complexity. Using models as in {\bf Step 1} we can  prove that such a surgery is admissible. Thus in this case $\hat{D}_M=0$.

\medskip

\noindent {\bf Case 2} If all $r_j$'s are  $=1$ then $H$ is a cyclic group, $M$ is a  homology space, so that $M\in \X$.

Suppose $r_1>1$.  We set $p:=p_1$ and
\[
G_p=\bigoplus_{j=1}^n A_p^{s_j}(q_j),\;\; 0<s_1\leq s_2\leq \cdots \leq s_k.
\]
We conclude as in {\bf Step 1} that all the  components of $\chi$   determined by the above decomposition of $G_p$ are nontrivial. By performing an $A_p$ surgery  on a mildly bad knot we obtain a manifold $M'$  satisfying
\[
\nu_{M'}=\nu_{M},\;\;r_1'=r_1-1,\;\;r_j'=r_j,\;\;\forall j=2,\cdots n.
\]
Thus $\kappa(M')< \kappa(M)$ and we conclude by induction.   We can now conclude that any  Dehn surgery which transforms an odd order ${\bQ}HS$ to an odd order ${\bQ}HS$ is admissible. We define the complexity of    an  odd order ${\bQ}HS$ to be the product of the   complexities of the $p$-groups  it decomposes into.

\bigskip

\noindent {\bf Step 3} If ${\gq}_M=A_2^n(q)\oplus {\gq}_1$, where  ${\gq}_1$ is a linking form on a group of odd order, then $M\in \X$.  Denote by $\R_2'$ the family of such rational homology spheres. Define the complexity of such a manifold to be
\[
\kappa_M =2^{n+1}\kappa({\gq}_1).
\]
The considerations in {\bf Step 2} lead to the following  complexity reduction trick.

\begin{lemma} If $K$ is a  knot is an odd order ${\bQ}HS$  such that $K^\perp$ is a nontrivial subgroup, then there exists $K'\in K^\perp$ and a surgery on $K'$ producing an odd order ${\bQ}HS$  of smaller complexity. Moreover, such a surgery is admissible.
\label{lemma: reduce0}
\end{lemma}

 Suppose ${\gq}_M=A_2^n(q)\oplus {\gq}_1$ and $\chi$ is a nontrivial character of ${\gq}_M$. Then $\chi^\perp\neq 0$. Decompose
 \[
 \chi=\chi_0\oplus \chi_1,\;\;\chi_0\in A_2^n(q),\;\;\chi_1\in {\gq}_1
 \]
 It follows that  $\chi_1^\perp$ is a nontrivial subgroup in ${\gq}_1$.  Perform a  complexity reduction  surgery on a knot $K\in \chi_1^\perp\subset {\gq}_1$ as in Lemma \ref{lemma: reduce0}  to    conclude as we have done before that $\widehat{D}_M(\chi)=0$.

 \bigskip

\noindent {\bf Step 4}  If ${\gq}_M= \bigoplus_{k=1}^mA_2^{n_k}(q_k)\oplus {\gq}_1$, $n_1\geq n_2\geq \cdots \geq n_m>0$,  where  ${\gq}_1$ is a linking form on a group of odd order, then $M\in \X$. Denote by $\R_2$ the family of such rational homology spheres.  Define the complexity of such a manifold to be
\[
\hat{\kappa}_M:= \kappa\bigl(\bigoplus_{k=1}^mA_2^{n_k}(q_k)\bigr)=2^{n_1+\cdots +n_m +m}.
\]
For  every $M\in \R_2$ we define its model  $\tilde{M}$  to be a connected sum of lens spaces with the same linking form  as $M$. Again we will carry a proof  by induction on the complexity.  The basic complexity reduction technique is contained in the following lemma whose proof is deferred to the Appendix.

\begin{lemma}(a) Suppose  $ c\in A_2^s(q_1)\oplus A_2^r(q_2)$, $s\geq r>0$, . Then there exists $K\in c^\perp$ of the form
\begin{equation}
K= K_1\oplus K_2
\label{eq: good}
\end{equation}
where $K_2$ is a generator of $A_2^r(q_2)$.

\noindent (b) Suppose $M$ is  a rational homology sphere  such that ${\gq}_M= A_2^s(q_1)\oplus A_2^r(q_2)$ and $K$ is a knot in $M$ whose homology class satisfies (\ref{eq: good}). Then   there exists  a Dehn surgery on $M$    such that the resulting manifold $M'$ is  in  $\R_2'$ and has smaller complexity. More precisely, we can arrange so that
\[
{\gq}_{M'}= A_2^t(q)\oplus {\gq}_1
\]
where $t\leq r+s$ and ${\gq}_1$ is the linking form of some  odd order lens space.
\label{lemma: reduce1}
\end{lemma}

Let  $M\in \R_2$. Then we can write
\[
{\gq}_M= {\gq}_0\oplus {\gq}_1 :=\bigl(\bigoplus_{k=1}^mA_2^{n_k}(q_k)\oplus {\gq}_1.
\]
If $m=1$ then $M\in \X$ according to {\bf Step 3}. We can assume $m>1$.

Any nontrivial  character $\chi\in {\gq}_M$ decomposes  as
\[
\chi =\chi_0 +\chi_1,\;\;\chi_i\in {\gq}_i,\;\;i=0,1.
\]
Pick $K\in A_2^{n_1}(q_1)\oplus A_2^{n_2}(q_2)$ orthogonal to
$\chi_0$, and satisfying  (\ref{eq: good}). We want to perform a
complexity reduction surgery as in Lemma \ref{lemma: reduce1} but
we first must show that any such surgery is admissible.      This
can be seen    by performing this surgery on the model $\tilde{M}\in \X$. It produces a manifold of smaller complexity which by induction we know it is in $\X$, and thus proving that the surgery is admissible. We can now conclude as many times before that  $\widehat{D}_M(\chi)=0$. This shows that ${\gR}_2\subset \X$.

\bigskip

\noindent {\bf Step 5} {\bf Conclusion} Suppose $M$ is an arbitrary  ${\bQ}HS$. Then ${\gq}_M={\gq}_0\oplus {\gq}_1$, where ${\gq}_0$ is a  linking form on a $2$ group, and ${\gq}_1$ is a linking form on an odd order group.   The results in \cite[Theorem 0.1]{KK}  show that     if we add sufficiently many $A$'s to ${\gq}_0$ we obtain a linking form isomorphic to a direct sum  of $A$'s. Topologically  this means that we can find a connected  sum  $X$ of lens spaces of order $2^s$ such that $M\# X\in\R_2\subset \X$.      Thus $\widehat{D}_{M\#X}=0$. Since $\widehat{D}$ is additive  with respect to connected sums we deduce $\widehat{D}_M=0$. This concludes the proof of  Theorem \ref{th: main}. $\Box$.

\section{Final comments}
\setcounter{equation}{0}
The invariant introduced by Ozsv\'{a}th  and Szab\'{o} in \cite{OS}  satisfies the same surgery formula as the   modified  Seiberg-Witten invariant, and detects in the same fashion  the Casson-Walker invariant.    This shows that the strategy presented in this paper also answers a question in \cite{OS}. More precisely, their invariant  is equivalent to the modified Reidemeister torsion.

 If we consider the mod ${\bZ}$ reduction  of the  modified Seiberg-Witten invariant we deduce that
 \[
 \ssw^0_M(\si)= \frac{1}{8}KS_M(\si) \mod {\bZ},
 \]
where $KS_M(\si)$ denotes the  Kreck-Stolz invariant. It general it depends on the metric but its  mod $8{\bZ}$ reduction is metric  independent. Fix a $spin$ structure $\epsilon$. This choice allows us to think of $\t$ and $\sw$ as functions $H\ra {\bQ}$, $H:=H_1(M, {\bZ})$.

Denote by $\f_M$ the space of functions $f: H\ra {\bQ}/{\bZ}$. For each $h\in H$ define the finite difference operator
\[
\Delta_h: \f\ra \f,\;\;(\Delta_h f)(\si):=f(h\cdot \si)-f(\si).
\]
In \cite{Tu5} its is shown that  for every $h_1,h_2\in H$ we have
\[
{\bf lk}_M(h_1,h_2)=\Delta_{h_1}\Delta_{h_2} \t\mod {\bZ} .
\]
Since the constant functions are killed by $\Delta_\bullet$ we deduce
\[
\Delta_\bullet \t=\Delta_\bullet \t^0.
\]
Our main result now implies the following equality.
\[
 {\bf lk}_M(h_1,h_2)=-\frac{1}{8}\Delta_{h_1}\Delta_{h_2} KS_M\mod {\bZ},\;\;\forall h_1,h_2\in {\bZ}.
\]
It would be interesting to investigate whether this  identity has a higher dimensional counterpart.

In \cite{N3}  we  associated to each $spin$ structure $\epsilon$ on a rational homology sphere an invariant $c(\epsilon)\in {\bQ}/{\bZ}$ which was powerful enough to distinguish many lens spaces. We can now  identify it.  We have
\[
c(\epsilon)=\frac{1}{8}KS_M(\epsilon)\mod {\bZ}.
\]

\appendix

\section{Proofs of  some technical results}
\label{s: a}
\noindent {\bf Proof of Lemma \ref{lemma: c2}.}\hspace{.3cm}   A simple model of  $A_p$ surgery is the manifold  given by the   surgery diagram
\[
\mathcal{D}_n:=\Bigl\{ (K_1,p^s/q_1), (K_2, p^r/q_2), (K,  n)\Bigr\}, \;\;s\geq r>0
\]
where $K_1$ and $K_2$ are unlinked, unknots, and $K$ is a knot such that ${\bf Lk}(K,K_i)=\ell_i$, $-p^s/2 < \ell_1 < p^s/2$, $-p^r/2<\ell_2<p^r/2$. The linking matrix of this diagram is
\[
\left[\begin{array}{ccc}
\frac{p^s}{q_1} & 0 & \ell_1 \\
& &\\
0 &\frac{p^r}{q_2} &\ell_2\\
& &\\
\ell_1 & \ell_2 & n
\end{array}
\right].
\]
We can view $K$ as a knot in the connected sum of lens spaces $L(p^s,-q_1)\# L(p^r,-q_2)$. $K$ is a bad knot if and only if
\begin{equation*}
\frac{q_1\ell_1^2}{p^s}+\frac{q_2\ell_2^2}{p^r}=k\in {\bZ}.
\tag{$\ast$}
\label{tag: ast}
\end{equation*}
 The  knot is mildly bad if and only if $(p,\ell_2)=1$. Denote by $H$ the first homology group of the  $3$-manifold obtained by performing the surgery indicated by $\mathcal{D}_n$. The matrix
\[
B_n:=\left[\begin{array}{ccc}
p^s & 0 & q_1\ell_1 \\
& &\\
0 &p^r &q_2\ell_2\\
& &\\
\ell_1 & \ell_2 & n
\end{array}
\right]
\]
is a presentation matrix  for $H$ and has determinant
\[
\det (B_n)=p^{s+r}n-p^sq_2\ell^2-p^rq_1\ell_1^2= p^{s+r}\Biggl(n-\Bigl(\frac{q_1\ell_1^2}{p^s}+\frac{q_2\ell_2^2}{p^r}\Bigr)\Biggr)=p^{s+r}(n-k).
\]
Observe that  $|\det B_{n}|= p^{s+r}$ when $n=k\pm 1$.  Let $n=k+1$.

Rewrite the  condition (\ref{tag: ast})  as
\[
q_1\ell_1^2+p^{s-r}q_2\ell_2^2=p^sk
\]
Since $(q_1\ell_1,p)=1$ we deduce that $(q_2\ell_2,p)=1$. To   find $H$ we need to find the elementary divisors $d_1|d_2|d_3$ of  $B_{k+1}$. Clearly $d_1=1$. By looking at the  $2\times 2$ minor in the top left hand corner we deduce that $d_2| p^{s+r}$. On the other hand, if we look   the $2\times 2$ minor
\[
\left|\begin{array}{cc}
0 & q_2\ell_2\\
\ell_1 &  k+1
\end{array}\right|
\]
we deduce  that it is not divisible by $p$. Thus $d_2=1$ which shows that $H$ is a cyclic $p$-group of order $p^{s+r}$. $\Box$

\bigskip

\begin{figure}[ht]
\begin{center}
\epsfig{figure=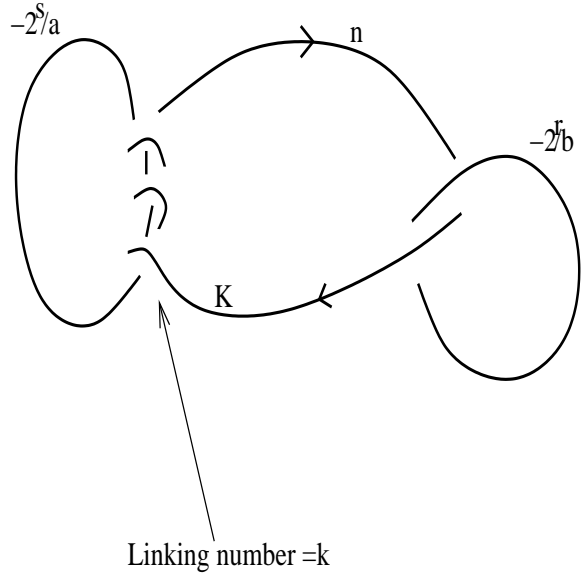, height=3in, width=3in}
\end{center}
\caption{\sl {Modeling a complexity reducing surgery}}
\label{fig: 6}
\end{figure}

\noindent {\bf Proof of Lemma \ref{lemma: reduce1}.} Part (a) is
elementary and  is left to the reader. For part (b) it suffices to
look   at a concrete realization of  the given homological data.
Any   homology class $K\in A_2^s(q_1)\oplus A_2^r(q_2)$ satisfying
(\ref{eq: good}) can be realized as a knot in a connected  sum of
lens spaces  $X:=L(2^s, a)\#L(2^r, b)$. We present $X$ as two unlinked unknots $K_1, K_2$ with surgery coefficients $-2^s/a$, $-2^r/b$, and $K$ as a knot  such that
\[
{\bf Lk}(K, {\bf g})=1,\;\;{\bf LK}(K,K_1)=k.
\]
Assume $K$  has an integral   surgery coefficient $n$ (see Figure \ref{fig: 6}).
 Slam-dunking $K_2$ over $K$ we obtain a surgery presentation with linking   matrix
\[
A:= \left[
\begin{array}{cc}
-2^s/a & k\\
k  & n+b/2^r
\end{array}
\right].
\]
The   first homology group $H$ of the manifold  described by this  surgery diagram admits  the presentation matrix
\[
B:=\left[
\begin{array}{cc}
-2^s & ak\\
2^rk  & 2^rn+b
\end{array}
\right].
\]
The order of this group is $|2^{r+s}n +2^sb+2^rak^2|$.  Pick $n$  to be any number such that $2^{r+s+1}$ does not divide the order of this group.  Then $H$ is a  cyclic group of the form ${\bZ}/2^t\oplus {\bZ}/(2m+1)$, $t\leq r+s$.  Its $\hat{\kappa}$-complexity is smaller than that of $A_2^s(q_1)\oplus A_2^r(q_2)$ $\Box$

\end{document}